\documentclass{article}

\usepackage{natbib, color}
\usepackage{graphicx, float}
\usepackage{multicol}
\usepackage{amsmath, fullpage, booktabs, amsthm}
\usepackage[colorlinks=TRUE, citecolor=blue]{hyperref}
\usepackage{amsfonts, bm}
\usepackage[margin=1in]{geometry}
\usepackage{setspace}

\usepackage[tablesfirst, nolists]{endfloat}

\allowdisplaybreaks

\newtheorem{Example}{Example}

\newcommand{\F}{{\mathcal{F}}}

\newcommand{\PP}{{\mathbf{P}}}
\newcommand{\EE}{{\mathbf{E}}}

\title{Martingale Approach to Gambler's Ruin Problem for Correlated Random Walks}

\author{Vladimir Pozdnyakov \thanks{Department of Statistics, University of Connecticut, 215 Glenbrook Road, Storrs, CT 06269-4120}\thanks{Corresponding author, vladimir.pozdnyakov@uconn.edu}}

\date{}

\begin{document}

\maketitle

\doublespacing


{
\noindent{\textbf{Abstract}}
The gambler's ruin problem for correlated random walks (CRW), both with and without delays,
is addressed using the Optional Stopping Theorem for martingales. We derive closed-form
expressions for the ruin probabilities and the expected game duration for CRW with
increments $\{1,-1\}$ and for symmetric CRW with increments $\{1,0,-1\}$ (CRW with delays).
Additionally, a martingale technique is developed for general CRW with delays.
The gambler's ruin probability for a game involving bets on two arbitrary patterns
is also examined.

\medskip
\noindent{\textbf{Keywords:}} Gambler's Ruin, Martingale, Optional Stopping Theorem, Correlated Random Walk

\medskip
\noindent{\textbf{2020 Mathematics Subject Classification:}} 60G42, 60J10
}

\section{Introduction}

A correlated random walk is a sequence of partial sums of random variables whose increments
form a Markov chain. Thus, it is a natural generalization of the standard random walk. The CRW is
also a special case of a Markov random walk (for example, \citet{Grama_etal2018} and \citet{Fuh2004}).
As it was noticed by \citet{GutStadtmüller2021}, a symmetric CRW with increments $\{1,-1\}$ can
also be viewed a  special case of an elephant random walk when
the elephant only remembers the very last step. In this paper we investigate a one-dimensional CRW;
a $d$-dimensional CRW is called the Gillis-Domb-Fisher  correlated random walk
(see \citet{ChenRenshaw1992} and \citet{Gillis1955}).

CRW has numerous applications in statistical physics, finance, ecology and other areas of research.
Here are just a few examples.  \citet{Allaart2004} employs CRW to  maximize the expected return
of a commodity seller. \citet{Gardner_etal2022} uses it to model short-term directional persistence
in animal movements. Applications to quantum walks in connection with quantum computing are used
as motivation  in \citet{Konno2009}. Our interest in random walks with delays
(that is, with $\{1,0,-1\}$ increments) is prompted by applications to sequential detection
in sensor networks with communication constraints (see \citet{Guerriero_etal2010}). In this work,
we consider a scenario where a sensor takes a measurement and, due to battery limitations,
proceeds as follows: if the measurement is too high, above a certain threshold,
it sends a 1 to the fusion center; if it is too low, it sends a -1.
When the measurement is within an acceptable range, the sensor does nothing,
which is interpreted as a passive response, or 0.

The gambler's ruin problem for CRW was first investigated in \citet{Mohan1955} with
the help of first-step analysis. The case of a two-state $k$-order Markov chain
was studied in \citet{MukherjeaSteele1987}  by the difference equations method.
However, as they noted, handling a large number  of difference equations is
difficult even for a small $k$. A numerical method  (via $n$-step transition probabilities)
for finding the ruin probabilities and the expected duration of the game is given by \citet{LalBhat1989}.
The gambler's ruin problem with delays for the independent increment case was solved in \citet{Gut2013}
with the help of martingales.

The paper proceeds as follows. In Section~\ref{two*patterns*section} we describe a game
that is used as an illustration of a Markovian gambler's ruin problem with delays.
In Section~\ref{crw*section} we derive closed-form  expressions for both the ruin probabilities
and the expected time until ruin for the correlated  random walk with increments $\{1, -1\}$.
Section~\ref{crw*delays*section} includes closed-form expressions for the ruin probabilities and
the expected time until ruin for symmetric CRW with  delays and a martingale treatment of the general case.
Section~\ref{two*pattern*ruin*section} provides a solution for the gambler's ruin problem in the case
of the two arbitrary pattern game. In the final section we briefly discuss some open questions.

\section{A Two-Pattern Problem}\label{two*patterns*section}

Let us introduce a simple game that naturally leads us to the gambler's ruin problem for correlated random walks with delays.
This game will be used to demonstrate the flexibility and limitations of our methodology.
The complete solution to this gambler's ruin problem will be provided in Section~\ref{two*pattern*ruin*section},
along with a solution for the general case.

A casino flips a fair coin. A gambler gets \$1, if he sees HH (including overlapping occurrences,
that is, for HHH he will get \$2). He loses a dollar, if it is TH.  Here are some facts that we know
about this game:
\begin{enumerate}
\item The expected number of occurrences of HH and TH in $n$ trials is $(n-1)/4$ for both patterns.
\item The expected waiting time until the first HH is 6.
\item The expected waiting time until the first TH is 4.
\item The expected waiting time until the first occurrence of either pattern is 3.
\item The probability that HH will occur before TH is .25.
\end{enumerate}

The first statement is a trivial fact that one can get by taking the expectation of a sum of appropriate
indicators. The other statements are more interesting. There is a long list of articles  on this topic,
the occurrence of patterns. An approach based on generating functions can be found in \citet{Feller1968}.
In a private communication, Donald Knuth indicated to us that the key observation---that
the expected waiting time until the first  occurrence of a pattern depends on self-overlapping of the pattern with
itself---was first made in \citet{Solovev1966}. A combinatorial method was
developed in \citet{GuibasOdlyzko1981}. More relevant to our approach is a nice martingale
technique introduced in \citet{Li1980}.

Now, let us ask a different question. Let $A$, $B$ be positive integers. What is the probability of
getting to \$$A$ before reaching -\$$B$? This is a variation of the famous gambler's ruin problem for
random walks that one can find in classical texts such as \citet{Feller1968}.
However, here we have two important differences from the standard setup: (1) the payments are
not independent, and (2) the payments take three values (1, 0, and -1).

We must acknowledge that there is an element of luck in the sequence of payments in this particular
game forming a Markov chain with three states $\{1,0,-1\}$. Usually, this is not the case.
One can easily verify that if we play the same game with HH and TT, the payments do not form a
Markov chain. However, as we demonstrate later, finding the ruin probability for a game with
two arbitrary patterns is closely related to the gambler's ruin problem for the CRW.

\section{Correlated Random Walk}\label{crw*section}

Consider a two-state Markov Chain $\{X_n\}_{n\geq 1}$ with state set $\{1, -1\}$, transition matrix
$$
\begin{pmatrix}
  p & 1-p \\
  1-q & q \\
\end{pmatrix},
$$
and initial distribution $(p_0, 1-p_0)$, where $0<p,q<1$ and $0\leq p_0\leq 1$. In what follows,
we assume that the initial distribution is always a stationary one, that is, in this case,
$p_0=(1-q)/(2-p-q)$. This assumption helps to eliminate the extra variable, $p_0$, in our
formulas. But all our derivations can be easily adjusted to an arbitrary initial distribution.

Let $S_n=X_1+\cdots +X_n$, $n\geq 1$ denote a correlated random walk, and let
$\F_n=\sigma(X_1, \dots, X_n)$, the $\sigma$-field generated by $\{X_n\}_{n\geq 1}$. For any integers $A,B>0$,
let us consider the following stopping time (w.r.t to the filtration $\{\F_n\}_{n\geq 1}$):
$$\tau=\inf\{k: S_k=A \mbox{ or }S_k=-B\}.$$
One important technical observation is that $\EE(\tau)<\infty$. It is a well-known fact, so we will skip the proof.
Then the {\it ruin} probabilities are:  $\alpha=\PP(S_\tau=A)$ and $\beta=\PP(S_\tau=-B)$.

When $p=1-q$, $\{X_n\}_{n\geq 1}$ are independent identically distributed (i.i.d.) random
variables. The gambler's ruin problem can be solved with the help of various techniques.
The most popular ones are first-step analysis and the martingale approach. A review of
both methods can be found in \citet{Steele2001}. In particular, the following martingales
are useful. For a symmetric random walk ($p=q=1/2$), they are $S_n$ and $S_n^2-n$.
For asymmetric random walks, they are $(q/p)^{S_n}$ and $S_n-(p-q)n$.
Then by a standard application of the Optional Stopping Theorem (OST)
one can show that, for example, in the symmetric case $\alpha=B/(A+B)$ and $\EE(\tau)=AB$.

The general Markovian case was first treated in \citet{Mohan1955} with the help of the first-step analysis.
To the best of our knowledge, nobody has used a martingale approach to address the gambler's ruin
problem for correlated random walks. So, our first task is to re-derive
the formulas from \citet{Mohan1955} with the help of the Optional Stopping Theorem. For that,
we need to construct an additive martingale (similar to $S_n-(p-q)n$ in the i.i.d. case).

The starting point is a formula for the conditional expectation $\EE(X_n|\F_{n-1})$. Observe first that
\begin{align*}
\EE(X_n|\F_{n-1})=&\EE(X_n|X_{n-1}=1)\mathbb{I}_{X_{n-1}=1}+\EE(X_n|X_{n-1}=-1)\mathbb{I}_{X_{n-1}=-1}\\
                 =&\big[1\cdot p+(-1)\cdot(1-p)\big]\frac{X_{n-1}+1}{2}+\big[1\cdot (1-q)+(-1)\cdot q\big]\frac{1-X_{n-1}}{2}\\
                 =&(p+q-1)X_{n-1}+(p-q).
\end{align*}
Next, if we define $Y_1=X_1$ and $Y_n=X_n-(p+q-1)X_{n-1}-(p-q)$ for $n\geq 2$, then
\begin{align}\label{additive*martingale}
M_n=&Y_1+\cdots+Y_n\nonumber\\
   =&(2-p-q)S_{n-1}+X_n-(p-q)(n-1)\\
   =&(2-p-q)S_{n}+(p+q-1)X_n-(p-q)(n-1)\nonumber
\end{align}
is a martingale w.r.t. the filtration $\{\F_n\}_{n\geq 1}$ as a sum of martingale-differences.

\subsection{Ruin Probability for Symmetric CRW} If $p=q$, then $$M_n=(2-2p)S_n+(2p-1)X_n.$$
The Optional Stopping Theorem immediately gives us that
$$0=\EE(M_1)=\EE (M_\tau)=(2-2p)\EE(S_\tau)+(2p-1)\EE(X_\tau).$$
Since $X_{\tau}=1$ on the set $\{S_{\tau}=A\}$ , and $X_{\tau}=-1$ on $\{S_{\tau}=-B\}$ we have
$$0=(2-2p)(\alpha A-\beta B)+ (2p-1)(\alpha-\beta).$$
Since we also have $\alpha+\beta=1$, we immediately get
\begin{equation}\label{alpha*symmetric}
\alpha=\frac{B-1+\frac{1}{2}\frac{1}{1-p}}{A+B-2+\frac{1}{1-p}}.
\end{equation}
Note that we have intuitive answers for special cases:
\begin{enumerate}
  \item if $A=B$ (symmetry), then $\alpha=1/2$,
  \item if $p=1/2$ (i.i.d.), then $\alpha=B/(A+B)$,
  \item if $p\to 1$, then $\alpha\to 1/2$ (everything is decided by the first flip),
  \item if $B\to\infty$, then $\alpha\to 1$.
\end{enumerate}

As mentioned above, adjusting to an arbitrary initial distribution $(\pi_1,1-\pi_1)$ is straightforward.
The only change is the value of $\EE(M_1)$. In this case, $\EE(M_1)=2\pi_1-1$. Consequently, we have
\begin{equation}\label{alpha*symmetric*with*initial}
\alpha=\frac{B-1+\frac{\pi_1}{1-p}}{A+B+\frac{2p-1}{1-p}}.
\end{equation}

\subsection{Expected Waiting Time until Ruin for Symmetric CRW}
The i.i.d. case tells us that we should try to find a {\it quadratic} martingale.
This suggests the following candidate: $$M_n=S_n^2+aS_nX_n-bn.$$ If it is possible
to find such values of $a$ and $b$ so $M_n$ forms a martingale, then the OST gives us
$$(1+a)-b=\EE(M_1)=\EE(M_\tau)=A^2\alpha+B^2\beta+a(A\alpha+B\beta)-b\EE(\tau).$$
Since we already know $\alpha$ and $\beta$, we have a formula for the expected waiting time
\begin{equation}\label{expectation*symmetric}
\EE(\tau)=\frac{1}{b}\big[b-(1+a)+A^2\alpha+B^2\beta+a(A\alpha+B\beta)\big].
\end{equation}
Now, using $\EE(X_n|\F_{n-1})=(2p-1)X_{n-1}$ and $X_n^2=1$ we get that
\begin{align*}
  \EE(M_n|\F_{n-1}) =& \EE(S_n^2+aS_nX_n-bn|\F_{n-1})\\
                    =& \EE((S_{n-1}+X_n)^2+a(S_{n-1}+X_n)X_n-bn|\F_{n-1})\\
                    =& S_{n-1}^2+\big[2(2p-1)+a(2p-1)\big]S_{n-1}X_{n-1}+\big[1+a-bn\big].
\end{align*}
Therefore, to make it a martingale we need $a$ and $b$ to satisfy the following system of equations:
$$
\left\{\begin{array}{l}
 2(2p-1)+a(2p-1) = a \\
   1+a-bn = -b(n-1)
\end{array}\right..
$$
This gives us $a=(2p-1)/(1-p)$ and $b=p/(1-p)$. Plugging these values into (\ref{expectation*symmetric}), after
some algebra, we obtain that
\begin{equation}\label{expectation*symmetric*final}
\EE(\tau)=AB\frac{1-p}{p}+(A+B)\frac{2p-1}{2p}.
\end{equation}
Here are some special cases:
\begin{enumerate}
  \item if $p=1/2$ (i.i.d.), then $\EE(\tau)=AB$,
  \item if $p=1$, then $\EE(\tau)=(A+B)/2$ (going straight to $A$ or $-B$ depending on initial value $X_1$),
  \item if $A=B$, then $\EE(\tau)=A^2\frac{1-p}{p}+A\frac{2p-1}{p}$.
\end{enumerate}

\subsection{Ruin Probability for Asymmetric CRW}
The i.i.d. case together with our first martingale construction suggests an {\it exponential} martingale
$$
M_n=\lambda^{S_n+aX_n},
$$
where $\lambda>0$ and $a\in \mathbb{R}$. Again, let us assume first that we are successful in finding
$\lambda>0$ and $a$ that make $M_n$ a martingale. Then in addition to $\alpha+\beta=1$, the OST
gives another linear equation:
\begin{equation}\label{alpha*asymmetric}
\EE(M_1)=\EE(M_\tau)=\lambda^{A+a}\alpha+\lambda^{-B-a}\beta.
\end{equation}

So, let us figure out these $\lambda$ and $a$. First, we have that
\begin{align*}
  \EE(M_n|\F_{n-1}) =& \EE(\lambda^{S_{n-1}+(1+a)X_n+ aX_{n-1}-aX_{n-1}}|\F_{n-1})
                    = M_{n-1}\EE(\lambda^{(1+a)X_n-aX_{n-1}}|\F_{n-1}).
\end{align*}
Can we make the conditional expectation in the last expression to be equal to 1? Yes. As before,
\begin{align*}
\EE(\lambda^{(1+a)X_n-aX_{n-1}}|\F_{n-1})
            =&\mathbin{\hphantom{+}}\big[ p\lambda^{(1+a)-a}+(1-p)\lambda^{-(1+a)-a}\big]\frac{X_{n-1}+1}{2}\\
             &+\big[(1-q)\lambda^{(1+a)+a}+q\lambda^{-(1+a)+a}\big]\frac{1-X_{n-1}}{2}\\
            =&\mathbin{\hphantom{+}}\big[ p\lambda +(1-p)\lambda^{-1-2a}-(1-q)\lambda^{1+2a}-q\lambda^{-1}\big]\frac{X_{n-1}}{2}\\
             &+\big[ p\lambda +(1-p)\lambda^{-1-2a}+(1-q)\lambda^{1+2a}+q\lambda^{-1}\big]\frac{1}{2}.
            \end{align*}
Next, equating the first coefficient to 0 and the second one to 2, we get the following non-linear system of
equations for $\lambda$ and $a$
$$
\left\{\begin{array}{l}
p\lambda+(1-p)\lambda^{-1-2a}=1 \\
(1-q)\lambda^{1+2a}+q\lambda^{-1}=1
\end{array}\right..
$$
Because of the i.i.d. martingale $(q/p)^{S_n}$, it is natural to try $\lambda=q/p$. And it works, $\lambda=q/p$ and
$$a=\frac{1}{2}\left[\log_{q/p}\frac{1-p}{1-q}-1\right]$$
is a solution. Note that $a$ itself is not really a part of equation~(\ref{alpha*asymmetric}), it enters as
$$\left(\frac{q}{p}\right)^a=\sqrt{\frac{(1-p)p}{(1-q)q}}.$$ Since our Markov chain starts from the stationary
distribution
\begin{align*}
\EE(M_1)=&\left(\frac{q}{p}\right)^{1+a}\frac{1-q}{2-p-q}+\left(\frac{q}{p}\right)^{-1-a}\frac{1-p}{2-p-q}
        =\frac{\sqrt{(1-q)(1-p)}}{2-p-q}\left[\sqrt{\frac{q}{p}}+\sqrt{\frac{p}{q}}\right].
\end{align*}
Finally, plugging it into (\ref{alpha*asymmetric}) and solving with respect $\alpha$ gives us
\begin{equation}\label{alpha*asymmetric*final}
\alpha=\frac{\left(\frac{q}{p}\right)^{-B}\sqrt{\frac{(1-q)q}{(1-p)p}}-
                   \frac{\sqrt{(1-q)(1-p)}}{2-p-q}\left[\sqrt{\frac{q}{p}}+\sqrt{\frac{p}{q}}\right]}
            {\left(\frac{q}{p}\right)^{-B}\sqrt{\frac{(1-q)q}{(1-p)p}}-\left(\frac{q}{p}\right)^{A}\sqrt{\frac{(1-p)p}{(1-q)q}}}.
\end{equation}
If $p=1-q$, $p\neq 1/2$ , then formula (\ref{alpha*asymmetric*final})
simplifies to a well-known expression for the asymmetric i.i.d. case (for example, \citep{Steele2001}):
$$
\alpha=\frac{\left(\frac{q}{p}\right)^{-B}-1}
            {\left(\frac{q}{p}\right)^{-B}-\left(\frac{q}{p}\right)^{A}}.
$$

\subsection{Expected Waiting Time until Ruin for Asymmetric CRW}
Everything is ready now. Applying the OST to the martingale~(\ref{additive*martingale}), we have
the following equation for $\EE(\tau)$:
$$\EE(M_1)=\EE(M_\tau)=(2-p-q)(A\alpha-B\beta)+(p+q-1)(\alpha-\beta)-(p-q)(\EE(\tau)-1).$$
Since we already know $\alpha$ (and $\beta=1-\alpha$) and $\EE(M_1)=\EE(X_1)=(p-q)/(2-p-q)$, the formula
for $\EE(\tau)$ immediately follows:
$$
\EE(\tau)=\frac{2-p-q}{p-q}(A\alpha-B\beta)+\frac{p+q-1}{p-q}(\alpha-\beta)+\frac{1-p-q}{2-p-q}.
$$

\section{Correlated Random Walk with Delays}\label{crw*delays*section}
As we mentioned in the introduction, accounting for the Markovian dependence is not enough to solve
our two-pattern problem. The payments are quite often equal to 0. Random walks associated with such
random variables are called random walks with delays. The gambler's ruin problem with delays for
the i.i.d. case can be solved  using standard techniques of martingales and the OST, see \citet{Gut2013}.
Here we extend these results to CRW with delays. We start with a symmetric CRW which can be viewed as
a natural generalization of the i.i.d. symmetric random walk with delays.

\subsection{Ruin Probability for Symmetric CRW with Delays}\label{ruin*symmetric*crw*section}
Consider a three-state Markov Chain $\{X_n\}_{n\geq 1}$ with state set $\{1, 0, -1\}$ and transition matrix
$$
\begin{pmatrix}
  p       & r  & 1-p-r \\
  (1-q)/2 & q  & (1-q)/2\\
  1-p-r   & r  & p \\
\end{pmatrix},
$$
where $0\leq p,q,r<1$ and $0\leq p+r\leq 1$. We assume that the initial distribution is stationary, that is, it is
given by
$$\left(\frac{1-q}{2(1+r-q)},\frac{r}{(1+r-q)},\frac{1-q}{2(1+r-q)}\right).$$ We will keep the same notation
for random walks, ruin probabilities, filtration and stopping time. We call this random walk symmetric, because
states 1 and -1 enter into the picture in a symmetric way, they are interchangeable. More specifically, for
states 1 and -1, the  probabilities of staying in the same state and the probabilities of switching to state 0
are the same. Additionally, the probabilities of switching from 0 to 1 and -1 are equal, too.

As before, the first step is to figure out the conditional expectation $\EE(X_{n}|\F_{n-1})$. We will use a
technique  similar to the one employed in the previous section. But note that in this case, $X_n^2$ is still a
random variable, and the expressions for indicators of events $\{X_n=k\}$, $k=1,0,-1$ now depend on
both $X_n$ and $X_n^2$. More specifically, we have
\begin{align*}
 &\mathbb{I}_{X_n=1}  =(X_n^2+X_n)/2, \quad\quad
 \mathbb{I}_{X_n=0}  =1-X_n^2,\quad\quad
 \mathbb{I}_{X_n=-1}  =(X_n^2-X_n)/2.
\end{align*}
As a consequence, we get
\begin{align*}
\EE(X_n|\F_{n-1})=&\EE(X_n|X_{n-1}=1)\mathbb{I}_{X_{n-1}=1}+\EE(X_n|X_{n-1}=0)\mathbb{I}_{X_{n-1}=0}+\EE(X_n|X_{n-1}=-1)\mathbb{I}_{X_{n-1}=-1}\\
                 =&\mathbin{\hphantom{+}}\big[p-(1-p-r)\big](X_{n-1}^2+X_{n-1})/2\\
                  &+\big[(1-q)/2-(1-q)/2\big](1-X_{n-1}^2)\\
                  &+\big[(1-p-r)-p\big](X_{n-1}^2-X_{n-1})/2\\
                 =&\mathbin{\hphantom{+}}(2p+r-1)X_{n-1}.
\end{align*}
Note that because of the symmetry of the transition matrix,  $X_{n-1}^2$ is not a part of the formula.

So, a natural candidate for another useful martingale is
$$M_n=S_n+aX_n=S_{n-1}+(a+1)X_n.$$ Using the expression for the conditional expectation, we
immediately get that $(a+1)(2p+r-1)=a$ or $a=(2p+r-1)/(2-2p-r)$.
Let us make two observations about the martingale. First, $q$ is not a part of the formula,
because in the symmetric case the delays do not affect the ruin probabilities. Second, note
that because of $p+r\leq1$ and $r<1$, we have $2-2p-r>0$, so no division by zero. The rest is standard.
The OST gives us
$$0=\EE(M_1)=\EE(M_\tau)= A\alpha-B\beta+\frac{2p+r-1}{2-2p-r}(\alpha-\beta).$$
Together with $\alpha+\beta=1$ it leads us to
\begin{equation}\label{alpha*symmetric*delays}
\alpha=\frac{B-1+\frac{1}{2-2p-r}}{A+B-2+\frac{2}{2-2p-r}}.
\end{equation}
As a check, let us take a look at some special cases:
\begin{enumerate}
  \item if $r=0$ (the state 0 is ignored), then we have formula (\ref{alpha*symmetric}),
  \item if $2p+r=1$ or $p=1-p-r$ (almost the i.i.d. case because of an additional symmetry), then $S_n$ itself a martingale, and $\alpha=B/(A+B)$
  \item if $A=B$ (symmetry), then $\alpha=1/2$.
\end{enumerate}
In fact, because the delays do not really affect the ruin probabilities, the gambler's ruin
problem for a symmetric CRW with delays is equivalent to the gambler's ruin
problem  for a two-state symmetric CRW with the probability of staying in a state that is equal to
$p+r/2$.
\subsection{Expected Waiting Time until Ruin for Symmetric CRW with Delays}
Finding $\EE(\tau)$ is more challenging, because the delays play a significant role.
At this point, we have some experience in building useful martingales, and this
experience suggests the following candidate:
$$M_n=S_n^2+aS_nX_n+bX_n^2-cn.$$
If we are successful in finding the right constants $a$, $b$, and $c$, the OST immediately
gives us a formula for the expected waiting time until ruin:
\begin{equation}\label{expectation*symmetric*delays}
\EE(\tau)=\frac{1}{c}\big[A^2\alpha+B^2\beta+a(A\alpha+B\beta)+b-\EE(M_1)\big].
\end{equation}

So, the first step is to figure out the conditional expectation $\EE(X_n^2|\F_n)$. Using the
expressions for indicators, we have that
\begin{align*}
\EE(X_n^2|\F_{n-1}) =&\mathbin{\hphantom{+}}\big[1-r\big](X_{n-1}^2+X_{n-1})/2\\
                  &+\big[1-q\big](1-X_{n-1}^2)\\
                  &+\big[1-r\big](X_{n-1}^2+X_{n-1})/2\\
                  =&\mathbin{\hphantom{+}}(1-q)+(q-r)X_{n-1}^2.
\end{align*}
Next, we get that
\begin{align*}
\EE(M_n|\F_{n-1}) =& \EE((S_{n-1}+X_n)^2+a(S_{n-1}+X_n)X_n+bX_n^2-cn|\F_{n-1})\\
                  =& \EE(S_{n-1}^2+2S_{n-1}X_n+X_n^2 +aS_{n-1}X_n+aX_n^2+bX_n^2-cn|\F_{n-1})\\
                  =&S_{n-1}^2+(2+a)(2p+r-1)S_{n-1}X_{n-1}\\
                   &\quad\quad +(1+a+b)(q-r)X_{n-1}^2\\
                   &\quad\quad +(1+a+b)(1-q)-cn.
\end{align*}
To make it a martingale $a$, $b$ and $c$ have to satisfy the following linear system:
$$
\left\{\begin{array}{l}
(2+a)(2p+r-1)=a \\
(1+a+b)(q-r)=b \\
(1+a+b)(1-q)=c
\end{array}\right..
$$
That leads us to
$$ a=\frac{2(2p+r-1)}{2-2p-r}, \quad\quad b=\frac{(q-r)(2p+r)}{(1+r-q)(2-2p-r)}, \quad\quad
c=\frac{(1-q)(2p+r)}{(1+r-q)(2-2p-r)}.
$$
Note that because of the restrictions on $p$, $q$, and $r$, there is no division by zero. The last step
is to calculate $\EE(M_1)$. Since we start from the stationary distribution, $\EE(X_1^2)=(1-q)/(1+r-q)$,
and, as a result, we get
\begin{align*}
\EE(M_1) =&\EE(S_{1}^2+aS_{1}X_1+bX_1^2-c)\\
         =&\EE(X_{1}^2+aX_{1}X_1+bX_1^2-c)\\
         =&(1+a+b)\frac{1-q}{1+r-q}-c.
\end{align*}
Finally, formula (\ref{expectation*symmetric*delays}) leads us to
\begin{equation}\label{expectation*symmetric*delays*final}
\EE(\tau)=AB\frac{(1+r-q)(2-2p-r)}{(1-q)(2p+r)}+(A+B)\frac{(2p+r-1)(1+r-q)}{(1-q)(2p+r)}+\frac{(q-r)r}{(1-q)(1+r-q)}.
\end{equation}
As always, let us consider special cases:
\begin{enumerate}
  \item if $r=0$ (the state 0 is ignored), then $q$ will be canceled out, and we get formula (\ref{expectation*symmetric*final}),
  \item if $2p+r=1$ or $p=1-p-r$ (balanced switches between 1 and -1), then
  $$\EE(\tau)=AB\frac{1+r-q}{1-q}+\frac{(q-r)r}{(1-q)(1+r-q)},$$
  \item if $2p+r=1$  and $r=q$ (the i.i.d. symmetric random walk with delays), then $\EE(\tau)=AB/(1-q)$, a formula from \citet{Gut2013}.
\end{enumerate}

\subsection{Ruin Probability for General CRW with Delays}
Finding closed-form formulas in a general case might be a challenging task. The transition matrix  has 6 free parameters.
However, it is possible to develop a method that will work, at least numerically, for almost any
Markov chain.

In this section  $\{X_n\}_{n\geq 1}$ is a three-state Markov Chain with state set $\{1, 0, -1\}$ and transition matrix
$$
\begin{pmatrix}
  p       & p_0  & 1-p-p_0 \\
  r_1 & r  & 1-r-r_1\\
  1-q-q_0   & q_0  & q \\
\end{pmatrix},
$$
where $0\leq p,p_0,q,q_0,r,r_1\leq 1$ and $0\leq p+p_0,r+r_1,q+q_0\leq 1$. We will drop here the assumption that the initial
distribution is stationary. The initial distribution is given by
$(\pi_1,\pi_0,1-\pi_0-\pi_1)$,
where $0\leq\pi_1,\pi_0\leq 1$ and $0\leq \pi_1+\pi_0\leq 1$. We will employ the same notation for random walks,
ruin probabilities, filtration and stopping time.

The candidate for a martingale in the general case, as before, takes an exponential form:
$$M_n=\lambda^{S_n+aX_n^2+bX_n}.$$
If one can find $\lambda>0$, $\lambda\neq 1$, $a$ and $b$ that make $M_n$ a martingale, then the OST
gives us
$$\EE (M_1)=\EE (M_\tau)=\lambda^{A+a+b}\alpha+\lambda^{-B+a-b}(1-\alpha).$$
Therefore, we have
\begin{equation}\label{alpha*asymetric*delays}
\alpha=\frac{\EE(M_1)-\lambda^{-B+a-b}}{\lambda^{A+a+b}-\lambda^{-B+a-b}}.
\end{equation}
So, let us try to find such  $\lambda$, $a$ and $b$. First, note that
\begin{align*}
  \EE(M_n|\F_{n-1}) =& \EE\left[\lambda^{S_{n-1}+aX_{n-1}^2+bX_{n-1}+a(X_n^2-X_{n-1}^2)+(b+1)X_n-bX_{n-1}}|\F_{n-1}\right] \\
                   =& M_{n-1}\EE\left[\lambda^{a(X_n^2-X_{n-1}^2)+(b+1)X_n-bX_{n-1}}|\F_{n-1}\right].
\end{align*}
Thus we need $\lambda$, $a$ and $b$ that make the conditional expectation of
$$Y_n=\lambda^{a(X_n^2-X_{n-1}^2)+(b+1)X_n-bX_{n-1}}$$
to be equal to 1. Now, our usual tricks lead us to
\begin{align*}
\EE(Y_n|\F_{n-1})=&\EE(Y_n|X_{n-1}=1)\mathbb{I}_{X_{n-1}=1}+\EE(Y_n|X_{n-1}=0)\mathbb{I}_{X_{n-1}=0}+\EE(Y_n|X_{n-1}=-1)\mathbb{I}_{X_{n-1}=-1}\\
=&\EE(Y_n|X_{n-1}=1)\frac{X_{n-1}^2+X_{n-1}}{2}+\EE(Y_n|X_{n-1}=-1)\frac{X_{n-1}^2-X_{n-1}}{2}\\
 &+\EE(Y_n|X_{n-1}=0)(1-X_{n-1}^2)\\
=&\mathbin{\hphantom{+}}\Big[\EE(Y_n|X_{n-1}=1)+\EE(Y_n|X_{n-1}=-1)-2\EE(Y_n|X_{n-1}=0)\Big]\frac{X_{n-1}^2}{2}\\
 &+\Big[\EE(Y_n|X_{n-1}=1)-\EE(Y_n|X_{n-1}=-1)\Big]\frac{X_{n-1}}{2}\\
 &+\Big[\EE(Y_n|X_{n-1}=0)\Big].
\end{align*}
Taking into account that
\begin{align*}
&\EE(Y_n|X_{n-1}=1)=p\lambda+p_0\lambda^{-a-b}+(1-p-p_0)\lambda^{-2b-1},\\
&\EE(Y_n|X_{n-1}=0)=r_1\lambda^{a+b+1}+r+(1-r-r_1)\lambda^{a-b-1},\\
&\EE(Y_n|X_{n-1}=-1)=(1-q-q_0)\lambda^{2b+1}+q_0\lambda^{-a+b}+q\lambda^{-1},
\end{align*}
we finally get the following system of non-linear equations
\begin{equation}\label{system*asymetric*delays}
\left\{\begin{array}{lr}
p\lambda+p_0\lambda^{-a-b}+(1-p-p_0)\lambda^{-2b-1}&=1\\
r_1\lambda^{a+b+1}+r+(1-r-r_1)\lambda^{a-b-1}& =1\\
(1-q-q_0)\lambda^{2b+1}+q_0\lambda^{-a+b}+q\lambda^{-1}&=1\\
\end{array}\right..
\end{equation}
\begin{Example}\label{Example}
{\rm Consider a three-state Markov chain with state set $\{1, 0, -1\}$, transition matrix
$$
\begin{pmatrix}
  1/2      & 1/4  & 1/4 \\
  1/3 & 1/3  & 1/3\\
  1/8   & 1/8  & 3/4 \\
\end{pmatrix},
$$
and the initial distribution $(0,1,0)$. Then one can show that $\lambda=13/10$,
$a=-\log_{13/10}\sqrt{81/65}$, and $b=\log_{13/10}\sqrt{20/13}$ solve system~(\ref{system*asymetric*delays}).
     Also in this case $\EE (M_1)=1$ . Plugging these values into formula~(\ref{alpha*asymetric*delays})
we get
$$ \alpha=\frac{(18/13)(13/10)^B-1}{(20/13)(13/10)^{A+B}-1}.
$$
As a simple check, set $A=B=1$, then obviously $\alpha=.5$. The formula gives the same answer.

One  practical recommendation is to do a re-parametrization,
and use $x=\lambda$, $y=\lambda^a$, and $z=\lambda^b$ instead of the original $\lambda$, $a$ and $b$.
Then {\tt Mathematica} easily finds solutions of the system.
}\end{Example}

This method might not work for some Markov chains, if their transition probability matrix exhibits  a certain symmetry.
Recall that in the i.i.d. case, $(q/p)^{S_n}$ is useless, if $p=q$.  In fact, our two-pattern
problem is a case of this nature.

\subsection{Expected Waiting Time until Ruin for General CRW with Delays}

Using the same technique as before one can show that
$$
\EE(X_n|\F_{n-1})=l_2X_{n-1}^2+l_1X_{n-1}+1_0,
$$
where
$$
l_2=p-q+(p_0-q_0)/2+1-r-2r_1,\quad\quad l_1=p+q+(p_0+q_0)/2-1,\quad\quad l_0=2r_1+r-1,
$$
and
$$
\EE(X_n^2|\F_{n-1})=k_2X_{n-1}^2+k_1X_{n-1}+k_0,
$$
where
$$
k_2=r-(p_0+q_0)/2,\quad\quad k_1=(q_0-p_0)/2,\quad\quad k_0=1-r.
$$
Remember that for the asymmetric i.i.d. case a simple martingale $S_n-cn$ works.
Therefore, it makes sense to try a {\it linear} martingale
$$
M_n=S_n+aX_n^2+bX_n-cn.
$$
Using our expressions for $\EE(X_n|\F_{n-1})$ and $\EE(X_n^2|\F_{n-1})$
we can show that if $a$, $b$, and $c$ satisfy
\begin{equation}\label{system*asymertic*delays*expectation}
\left\{\begin{array}{lr}
ak_2+(1+b)l_2&=a\\
ak_1+(1+b)l_1&=b\\
ak_0+(1+b)l_0&=c\\
\end{array}\right.,
\end{equation}
then $M_n$ is a  martingale. Moreover, if $c\neq 0$, then we get a formula the expected
time until ruin:
$$
\EE(\tau)=\frac{1}{c}\big[(A+a+b)\alpha+(-B+a-b)(1-\alpha)-\EE(M_1)\big],
$$
where
$$
\EE(M_1)=\pi_1(1+a+b)+(1-\pi_1-\pi_0)(-1+a-b)-c.
$$
If we consider the Markov chain in Example~\ref{Example}, then our method will lead
to the following martingale:
$$
M_n=S_n-\frac{2}{5}X_n^2+\frac{37}{45}X_n+\frac{4}{15}n.
$$

\section{Gambler's Ruin Probability for Two-Pattern Game}\label{two*pattern*ruin*section}

In this section, we first provide a complete solution of the HH vs TH gambler's ruin problem. Then,
we describe a martingale method that can be used to find the ruin probabilities for an arbitrary
two-pattern game. Finally, we will use the method to re-derive the ruin probability for the HH vs TH game.

\subsection{HH vs TH Gambler's Ruin Problem}
Let $\{\xi_n\}_{n\geq 1}$ be a sequence of independent Bernoulli random variables with a probability of success 1/2.
Then for $n\geq 2$ our payments are given by
$$
X_n=\left\{\begin{array}{rl}
             1,& \mbox{ if } \xi_{n-1}=1, \xi_n=1, \\
             0,& \mbox{ if }  \xi_n=0,\\
             -1,& \mbox{ if } \xi_{n-1}=0, \xi_n=1.
           \end{array}
\right.
$$
In the general case of two patterns of length 2 (for instance, HH and TT), the associated
sequence of payments typically does not form a Markov chain.
But for HH and TH, $\{X_n\}_{n\geq 2}$ is a first-order Markov chain with the transition matrix
$$
\begin{pmatrix}
 1/2      & 1/2  & 0 \\
 0 & 1/2  & 1/2\\
1/2  & 1/2  & 0 \\
\end{pmatrix},
$$
and the initial distribution $(1/4, 1/2, 1/4)$. This Markov chain is not symmetric in the sense
of Section~\ref{ruin*symmetric*crw*section}, but it does have a very special structure. As a result,
one can easily verify that for this transition matrix, there is no solution of
system~(\ref{system*asymetric*delays}) such that $\lambda\neq 1$.

However, we can construct a linear one. Note that here $S_n=X_2+\cdots+X_n$. Since
$$
\EE(X_n|\F_{n-1})=X_{n-1}^2-\frac{1}{2},
$$
and
$$
\EE(X_n^2|\F_{n-1})=\frac{1}{2},
$$
using our approach (or by solving linear system~(\ref{system*asymertic*delays*expectation})) one can
verify that $M_n=S_n+X_n^2,\quad n\geq 2$ is a martingale. Applying the OST to this martingale gives us
$$
\alpha=\frac{B-1/2}{A+B}.
$$
A simple check in this case is when $A=B=1$. The formula gives 1/4. It makes sense, because
the only way to get to 1 is to have HH right away.

Since the linear martingale does not contain an $n$ term, it cannot be used to compute $\EE(\tau)$.
But then we can try to find a quadratic martingale similar to the one used for the symmetric CRW with
delays. Indeed, one can show that
$$
M_n=S_n^2+2S_nX_n^2+2X_n^2-n/2
$$
is a martingale. Applying the OST to this martingale gives us that
$$
\EE(\tau)=2AB+B-A+1.
$$

\subsection{Gambler's Ruin Probability for Arbitrary Two-Pattern Game}

As mentioned above, for an arbitrary two-pattern game the sequence of payments typically is not a Markov chain.
For example, suppose that we play this game with HH and TT. As before, let $\{X_n\}_{n\geq 2}$ be the sequence of payments,
where we get 1 for HH, -1 for TT, and 0, otherwise. Then for any $n\geq 2$, we have
$$\PP(X_{2n+1}=1|X_{2n}=0,X_{2n-1}=0, \dots, X_3=0,X_2=1)=\frac{1}{2},$$
but
$$\PP(X_{2n+1}=1|X_{2n}=0,X_{2n-1}=0, \dots, X_3=0,X_2=-1)=0.$$

However, we can still find the ruin probability using CRW. For the rest of the section, let
$\{\xi_n\}_{n\geq 1}$ be an i.i.d. sequence of letters from a finite alphabet. A
sequences of letters is called a {\it pattern}. A pattern that has no letters is
called a {\it void pattern}. Consider two patterns
$$P=a_1a_2\dots a_k,\quad\quad\mbox{and}\quad\quad Q=b_1b_2\dots b_m,$$
where $k,m\geq 1$. Assume that $P$ and $Q$ are not connected subsequences of each other.
Our game is the same as before: we get \$1 for every (overlapping)
occurrence of $P$ and we lose a dollar if we see $Q$.

Let $\{X_n\}_{n\geq 1}$ be the sequence of only 1 or -1 payments (that is, 0s are ignored).
Then $\{X_n\}_{n\geq 1}$ is a two-state Markov chain whose initial distribution and
transition probabilities can be computed with the help of the theory developed for the occurrence
of patterns. More specifically, let
$$P^*=\left\{\begin{array}{rl}
             a_2\dots a_k,&  k\geq 2, \\
             \mbox{void pattern},& k=1,
           \end{array}
\right.
$$
and
$$Q^*=\left\{\begin{array}{rl}
             b_2\dots b_m,&  m\geq 2, \\
             \mbox{void pattern},& m=1.
           \end{array}
\right.
$$
Then the initial probability
$$\pi_1=\PP(\mbox{pattern $P$ will occur before $Q$ in the sequence$\{\xi_n\}_{n\geq 1}$}).$$
The transition probability of going from state 1 to state 1 is given by
$$p=\PP(\mbox{pattern $P$ will occur before $Q$ in the sequence $\{\xi_n\}_{n\geq 1}$ given that it starts from $P^*$}),$$
and the transition probability of going from state -1 to state -1 is given by
$$q=\PP(\mbox{pattern $Q$ will occur before $P$ in the sequence $\{\xi_n\}_{n\geq 1}$ given that it starts from $Q^*$}).$$
All three probabilities can be found using, for example, Theorem 3.1 in \citet{Li1980}.
To find each probability, it is required to solve a linear system with three equations.
Then, we simply apply the appropriate formula from Section~\ref{crw*section}.

To close the loop let us apply this method to  the HH vs TH game. As we mentioned in Section~\ref{two*patterns*section},
in this case $\pi_1=.25$. Additionally, here $P^*=Q^*=H$, and $p=q=1/2$. This means that the non-zero
payments in the HH vs. TH game form a sequence of independent symmetric random variables,
primarily because the coin is fair. Plugging these values
into formula~(\ref{alpha*symmetric*with*initial}) leads us again to
$$\alpha=\frac{B-1/2}{A+B}.$$

\section{Concluding Remarks} The expected time until the gambler's ruin in a game with two
arbitrary patterns remains an open question. The method based on the occurrence of
patterns ignores zero-payments, making it unclear how it can be applied to this problem.
Additionally, our technique relies on the fact that indicators of events like $X_n=1$ can
be expressed in terms of $X_n$. Consequently, it is likely unsuitable for Markov chains
with more than three states. The distribution of the waiting time until ruin for
CRW with delays is an open question as well. Finally, the gambler's ruin problem for elephant walks is also
unresolved. A comprehensive review of the results on elephant random walks can be found
in the recent work by \citet{GutStadtmuller2023}.


\bibliographystyle{mcap}
\bibliography{crwm}

\end{document}